\documentclass[12pt]{amsart}
\usepackage{amssymb}
{\catcode`\@=11\gdef\logo@{}}
\begin{document}
\baselineskip=16.5pt
\newcommand\irta{\mathcal A}
\newcommand\irtd{\mathcal D}
\newcommand\irtz{\mathcal Z}
\newcommand\irtp{\mathcal P}
\newcommand\irtb{\mathcal B}
\newcommand\irty{\mathcal Y}
\newcommand\irtu{\mathcal U}
\newcommand\irtl{\mathcal L}
\newcommand{\stick}%
{\mbox{{\hspace{0.4ex}$\mbox{\raisebox{1ex}{$\bullet$}}%
\hspace*{-0.94ex}|$\hspace{0.6ex}}}}

\centerline{{\bf STRONGLY ALMOST DISJOINT FAMILES, II}}
\bigskip
\centerline{by}
\bigskip
\centerline{A. Hajnal$^{1)}$, I. Juhasz$^{2)}$ and S. Shelah$^{3)}$}

\bigskip\bigskip\bigskip

{\it Abstract.}\quad The relations $M(\kappa,\lambda,\mu)\to
 B\,\,[\text{resp.\,} 
B(\sigma)]$ meaning that if $\mathcal A\subset[\kappa]^\lambda$ 
with $|\irta|=\kappa$ is $\mu$-almost disjoint then $\irta$ has 
property $B$ [resp. has a $\sigma$-transversal] had been 
introduced and studied under GCH in [EH]. Our two main results 
here say the following:

Assume GCH and $\varrho$ be any regular cardinal with a 
supercompact [resp. 2-huge] cardinal above $\varrho$. Then 
there is a $\varrho$-closed forcing $P$ such that, in 
$V^P$, we have both GCH and
$M(\varrho^{(+\varrho+1)},\varrho^+,\varrho)\nrightarrow
B\,\,[\text{resp.\,}M(\varrho^{(+\varrho+1)},\lambda,\varrho) \nrightarrow 
B(\varrho^+)$ for all $\lambda\le\varrho^{(+\varrho+1)}]$.

These show that, consistently, the results of [EH] are sharp. 
The necessity of using large cardinals follows from the results of [HJSh] 
and [BDJShSz].

\vskip4truecm



Research supported by\newline
1) NSF grant DMS-9704477;\newline
2) NSF grant DMS-9704477 and OTKA grant 25745;\newline
3)  NSF grant DMS-9704477 and the Israel Science Foundation 
founded by the Israel Academy of Sciences and Humanities. 
Publication no. 697.\newline
\newpage

\centerline{{\bf \S 1. Introduction}}

\smallskip

The aim of this paper is to show that, assuming the existence of certain 
large cardinals, the results of [EH] are sharp. Let us recall these results, 
first their terminology.

If $\mu\le\lambda\le\kappa$ and $\sigma$
are infinite cardinals then $M(\kappa,\lambda,\mu)\to B(\sigma)$ \hfill\break
$[\text{resp.\,}M(\kappa,\lambda,\mu)\to B]$ abbreviates 
the following statement: Whenever $\irta\subset[\kappa]^\lambda$ with 
$|\irta|=\kappa$ is $\mu$-almost disjoint (in short: $\mu$-a.d.) 
then $\irta$ has a $\sigma$-transversal [resp. $\irta$ has property $B$]. 
Here $\irta$ is $\mu$-a.d. means that the intersection of any two members 
of $\irta$ has size $<\mu$; a $\sigma$-transversal of $\irta$ is a set 
$T$ such that $0<|A\cap T|<\sigma$ holds for every $A\in\irta$; and 
$\irta$ has property $B$ if there is a set $T$ with $\emptyset\ne A
\cap T\ne A$ for all $A\in\irta$.

One of the main results of [EH] (see also [W, Chapter 1]) is as follows:

\bigskip

{\bf 1.1. Theorem.} (GCH)\quad If $\varrho$ is any regular cardinal 
then for any $\lambda\le\kappa\le\varrho^{(+\varrho)}$ we have
$$
M(\kappa,\lambda,\varrho)\to B(\varrho^+).
$$

\bigskip

The natural question whether the restriction $\kappa\le\varrho^{(+\varrho)}$ 
is essential here had also been raised in [EH], especially because 
the following was also proved there.

\bigskip

{\bf 1.2. Theorem.} (GCH)\quad If $\varrho$ is regular then for any 
$\lambda\le\kappa$
$$
M(\kappa,\lambda,\varrho)\to B(\varrho^{++})
$$
is valid. So if also $\lambda>\varrho^+$ then $M(\kappa,\lambda
,\varrho)\to B$.

\bigskip

Concerning the above question it was much later shown in [HJSh] and 
[BDJShSz] that the restriction $\kappa\le\varrho^{(+\varrho)}$ in 1.1 
can be omitted if some weak $\square$-like principles hold in addition 
to GCH, hence e.g. if $V=L$. On the other hand, it was also shown 
in [HJSh] that the existence of a supercompact cardinal implies the 
consistency of $M(\aleph_{\omega+1},\aleph_1,\aleph_0)
\nrightarrow B$, hence also of $M(\aleph_{\omega+1},\aleph_1,\aleph_0)
\nrightarrow B(\aleph_1)$, with GCH. The appearence of large cardinals 
here is of course essential because one has to negate the above 
mentioned $\square$-like principles.

Our first main result generalizes this negative result from $\varrho=
\aleph_0$ to any regular cardinal $\varrho$. This was not immediate because
the method of proof used in [HJSh] does not apply if $\varrho>\aleph_0$, 
so a new ingredient was needed. The general result can be formulated as
follows. 

\bigskip

{\bf 1.3. Theorem.}\qquad 
Assume that GCH holds, $\varrho$ is any regular cardinal and $\kappa$
is a supercompact cardinal with $\varrho<\kappa$. Then there is a 
$\varrho$-closed notion of forcing $P$ such that, in $V^P$, 
we have GCH and 
$$
M(\varrho^{(+\varrho+1)}, \varrho^+,\varrho)\nrightarrow B.
$$

\bigskip

(Note that since $P$ is $\varrho$-closed, no cardinals or cofinalities 
will be changed in $V^P$ up to $\varrho$.)

Of course, we trivially have here again that $M(\varrho^{(+\varrho+1)},
\varrho^+,\varrho)\nrightarrow B(\varrho^+)$ holds as well, but 
the relations $M(\varrho^{(+\varrho+1)},
\lambda,\varrho)\nrightarrow B(\varrho^+)$ are not excluded for 
$\varrho^+<\lambda\le\varrho^{(+\varrho+1)}$ . Our second main 
result, formulated below, takes care of these.

\bigskip

{\bf 1.4. Theorem.}\quad  
Assume GCH, $\varrho$ is regular and 
$\kappa$ is 2-huge with $\varrho<\kappa$. Then there is a 
$\varrho$-closed notion of forcing P such that, in 
$V^P$, we have GCH and 
$$
M(\varrho^{(+\varrho+1)},
\lambda,\varrho)\nrightarrow B(\varrho^+)
$$
for all $\lambda\le\varrho^{(+\varrho+1)}$.

\bigskip

These results indeed show that, modulo some large cardinals, the 
results of [EH] are best possible. However, the question 
of exactly what large cardinals are needed, in particular if the rather 
large step from the supercompact of 1.3 to the 2-huge of 1.4 is necessary, 
remains open.

\bigskip

\centerline{{\bf \S 2. The proof of 1.3}}
\smallskip

We start by recalling the following simple result from [HJSh]:

\bigskip

{\bf 2.1. Lemma.}\quad Let $S\subset\kappa$ be a stationary set 
such that $\diamondsuit(S)$ holds and 
$\{A_\alpha\colon\alpha\in S\}$ be a family of 
infinite sets with $A_\alpha\subset\alpha$ for each 
$\alpha\in S$. Then we can find sets 
$B_\alpha\subset A_\alpha$ with $|B_\alpha|=|A_\alpha|$ for all 
$\alpha\in S$ so that the family $\{B_\alpha\colon\alpha\in S\}$
does not have property $B$.

\bigskip

Let us now fix the regular cardinal $\varrho$, and to simplify notation 
let us denote $\varrho^{(+\varrho+1)}$ by $\widehat\varrho$. Also, 
given two regular cardinals $\lambda$ and $\kappa$ with $\lambda<\kappa$ we 
set 
$$
S_\lambda^\kappa=\{\alpha\in\kappa\colon\text{cf}(\alpha)=\lambda\}.
$$
Thus, by 2.1, $M(\widehat\varrho,\varrho^+,\varrho)\nrightarrow B$ is 
valid if we can find a stationary set 
$S\subset S_{\varrho^+}^{\widehat\varrho}$ satisfying 
$\diamondsuit(S)$ and a $\varrho$-a.d. family 
$\{A_\alpha\colon\alpha\in S\}$ such that $A_\alpha\in [\alpha]^{\varrho^+}$ 
for each $\alpha\in S$. Note that, as is well-known, GCH implies 
$\diamondsuit(S)$
whenever $S\subset S_\lambda^\kappa$ is stationary.

So far, everything has been done as in [HJSh] for the case 
$\varrho=\aleph_0$. It is the following theorem that allows us 
to get the result for an arbitrary regular cardinal $\varrho$.

\bigskip

{\bf 2.2. Theorem.}\quad Let $\mu$ be a singular cardinal of 
cofinality $\varrho$ and such that $\mu=\mu^{<\varrho}$. Then there is a 
partial order $Q=Q(\mu)$ with properties (i)-(v)
below:
\begin{itemize}
\item[(i)] $Q$ is $\varrho$-closed;  
\item[(ii)] $Q$ is $\mu^+$-CC; 
\item[(iii)] $|Q|\le 2^{\mu}$; 
\item[(iv)] in $V^Q$, $\mu$ is collapsed to $\varrho$, and 
$\varrho^+=\mu^+$;
\item[(v)] there is, 
in $V^Q$, a set $X\in[\varrho^+]^{\varrho^+}$ 
such that for any set $H\in\irtp(\mu^+)\cap V$ we have 
$|H\cap X|^{V^P}<\varrho$ if and only if $|H|^V<\mu$. 
\end{itemize}
\bigskip

{\it Proof.}\quad Let us put $Q=Q_1\times Q_2$, where $Q_1$ is the 
natural $\varrho$-closed partial order that forces a map of 
$\varrho$ onto $\mu$, i.e. $q\in Q_1$ iff $q$ maps some $\alpha\in\varrho$ 
into $\mu$, and extension is the partial ordering. Moreover,
$$
Q_2=[\mu^+]^{<\varrho}\times [\mu^+]^{<\mu}
$$
with the following ordering: for $\langle a,A\rangle$, 
$\langle a',A'\rangle\in Q_2$ we have 
$\langle a,A\rangle\le\langle a',A'\rangle$ iff 
$a\supset a'$, $A\supset A'$ and 
$A'\cap (a\setminus a')=\emptyset$ hold.

Clearly, both $Q_1$ and $Q_2$ are $\varrho$-closed, hence so is 
$Q$, i.e. (i) holds.

To show (ii), let us first note that from $\mu=\mu^{<\varrho}$ we have 
$|Q_1|=\mu$ and so it suffices to prove that $Q_2$ is $\mu^+$-CC. Thus 
let $\{\langle a_i,A_i\rangle\colon i\in\mu^+\}\subset Q_2$; 
clearly we may assume that 
$|a_i\cup A_i|<\lambda$ holds for a fixed regular cardinal 
$\lambda<\mu$ for all $i\in\mu^+$. Now, for every 
$\gamma\in S_\lambda^{\mu^+}$ the set $B_\gamma=
(a_\gamma\cup A_\gamma)\cap\gamma$ is bounded in $\gamma$, i.e. 
there is an $f(\gamma)<\gamma$ with $B_\gamma\subset f(\gamma)$. 
So by Fodor's theorem there is a stationary set $S\subset S_\lambda^{\mu^+}$ 
on which $f$ takes the constant value $\alpha$. Using $\mu^{<\varrho}=\mu$ 
we may also assume that $a_\gamma\cap\gamma=
a_\gamma\cap\alpha=c$ for all $\gamma\in S$.

Let us now pick $\gamma,\delta\in S$ such that both 
$\gamma<\delta$ and $a_\gamma\cup A_\gamma\subset\delta$, this is 
possible because each $a_\gamma\cup A_\gamma$ is bounded in $\mu^+$, and 
set $a=a_\gamma\cup a_\delta$, $A=A_\gamma\cup A_\delta$. 
Clearly, we have $\langle a,A\rangle\in Q_2$ and we next show that 
$\langle a,A\rangle$ extends both 
$\langle a_\gamma,A_\gamma\rangle$ and $\langle a_\delta, A_\delta\rangle$. 
Indeed, this follows because $a\setminus a_\gamma=a_\delta
\setminus c\subset\mu^+\setminus\delta$ and $A_\gamma
\subset\delta$ imply $A_\gamma\cap(a\setminus a_\gamma)=\emptyset$, 
moreover $a\setminus a_\delta=a_\gamma\setminus c\subset\delta\setminus
\gamma$ and $A_\delta\subset\alpha\cup(\mu^+\setminus\delta)$ imply 
$A_\delta\cap(a\setminus a_\delta)=\emptyset$.

(iii) follows easily because $|Q_1|=\mu$ and 
$|Q_2|=(\mu^+)^{<\varrho}(\mu^+)^{<\mu}\le(\mu^+)^\mu=2^\mu$.

(iv) is again trivial because $Q_1$ collapses $\mu$ to $\varrho$ and 
by (ii) $\mu^+$ is preserved.

Finaly, to see (v), let 
$G=G_1\times G_2$ be $Q$-generic over $V$ and set, in $V[G]$,
$$
X=\cup\{a\colon(\exists A)(\langle a,A\rangle\in G_2)\}.
$$   
Clearly, for every $\alpha\in\mu^+$ the set
$$
D_\alpha=
\{\langle a,A\rangle\in Q_2\colon a\setminus\alpha\ne\emptyset\}
$$
is dense in $Q_2$ and so $X$ is unbounded in $\mu^+=\varrho^+$, i.e. 
$X\in [\varrho^+]^{\varrho^+}$.

Now, if $H\in [\mu^+]^{<\mu}\cap V$
then again
$$
D_H=\{\langle a,A\rangle\in Q_2\colon H\subset A\}
$$
is dense in $Q_2$ because $\langle a,A\cup H\rangle\le\langle a,A\rangle$ 
for each $\langle a,A\rangle\in Q_2$. But then $G_2\cap D_H\ne\emptyset$, 
and if $\langle a,A\rangle\in G_2\cap D_H$ then we clearly have 
$X\cap H\subset X\cap A\subset a$, hence $|X\cap H|<\varrho$.

If, on the other hand, $H\subset \mu^+$, $H\in V$
and $|H|\ge\mu$ then clearly
$$
E_H=\{\langle a,A\rangle\in Q_2\colon a\cap H\ne\emptyset\}
$$
is dense in $Q_2$. Now, if we had $|X\cap H|<\varrho$
then by (i) we also had $X\cap H\in V$ and so $H\setminus X\in V$ and 
$|H\setminus X|\ge\mu$. This, however, contradicts the denseness of 
$E_{H\setminus X}$.

The following corollary is now immediate.

\bigskip

{\bf 2.3. Corollary.}\quad With the assumptions of 2.2, we can, in
$V^Q$, associate with every ground model set $A\in V$ with $|A|=
\varrho^+=\mu^+$ a subset $A^*\in [A]^{\varrho^+}$ such that for any
set $B\in V$ we have $|A^*\cap B|<\varrho$ iff $|A\cap B|^V<\mu$. 
In particular, if $\irta$ is a $\mu$-a.d. family of sets of size 
$\mu^+$ in $V$ then $\irta ^*=
\{A^*\colon A\in\irta\}$ is a $\varrho$-a.d. family of sets of size 
$\varrho^+$ in $V^Q$.

\bigskip

{\it Proof.}\quad Let $h\colon\mu^+\to A$ be a bijection of $\mu^+$ onto 
$A$ in $V$. Clearly,
$$
A^*=\{h(\xi)\colon\xi\in X\}
$$
is as required by (v) of 2.2.

Let us now return to the proof of 1.3. Let us put $\lambda=
\kappa^{(+\varrho)}$ and since $\kappa$ is $\lambda^+$-supercompact 
we can fix a normal, $\kappa$-complete ultrafilter $\irtu$ on 
$[\lambda^+]^{<\kappa}$. Using GCH we get $(\lambda^+)^\varrho=
\lambda^+$, hence we may also fix a bijection $G$ of $[\lambda^+]^\varrho$ 
onto $\lambda^+$. Standard reflection arguments and Solovay's Theorem 2 
from [S] then imply the existence of a set $A\in\irtu$ 
such that
\begin{itemize}
\item[(i)]\quad the map $P \mapsto \cup P$ is one-one on $A$;
\item[(ii)]\quad  each $P\in A$ is $G$-closed;
\item[(iii)]\quad $P\cap\kappa$ is an inaccessible cardinal and
$$
tp(P)=(P\cap\kappa)^{(+\varrho+1)}
$$
for each $P\in A$.
\end{itemize}
Now the set $S_1=\{\cup P\colon P\in A\}$ is clearly stationary 
in $\lambda^+$ since $\irtu$ is normal and, by (i),
we have $A=\{P_\alpha\colon\alpha\in S_1\}$ where $\cup P_\alpha=\alpha$ 
for $\alpha\in S_1$.

Let us now consider the map $\alpha\mapsto P_\alpha\cap\kappa$ on 
$S_1$. Then by (iii) we have a fixed inaccessible cardinal $\tau$ such 
that
$$S=\{\alpha\in S_1\colon P_\alpha\cap\kappa=\tau\}$$
is also stationary. We claim that the family 
$\{P_\alpha\colon\alpha\in S\}\subset[\lambda^+]^{\tau^{(+\varrho+1)}}$ is 
also $\tau^{(+\varrho)}$-a.d. Indeed, if $\alpha,\beta\in S$ are distinct 
and $|P_\alpha\cap P_\beta|\ge\tau^{(+\varrho)}$ held then by (ii) 
we also had $|P_\alpha\cap P_\beta|=\tau^{(\varrho+1)}$, using that 
$P_\alpha\cap P_\beta$ is $G$-closed. This, however contradicts that 
tp$(P_\alpha)=$tp$(P_\beta)=\tau^{(\varrho+1)}$ and $\cup P_\alpha=
\alpha\ne\cup P_\beta=\beta$.

Note that the singular cardinal $\mu=\tau^{(+\varrho)}$ satisfies 
the conditions of 2.2, hence in $V^{Q(\mu)}$ the family $\{P_\alpha^*\colon
\alpha\in S\}\subset[\lambda^+]^{\varrho^+}$ is $\varrho$-a.d., according to 
2.3. All that remains to be done is now to do a further $\varrho$-closed 
forcing that turns $\lambda^+$ into $\widehat\varrho$ and preserves both 
GCH and the stationarity of $S$. This job will clearly be done by e.g. 
$L v(\kappa,\varrho^{++})$, i.e. the Levy collapse of $\kappa$ to 
$\varrho^{++}$ in $V^{Q(\mu)}$. Then $P=Q(\mu)*L v(\kappa,\varrho^{++})$ 
is a $\varrho$-closed forcing such that $V^P$ satisfies GCH, 
moreover, in $V^P$,  
$\{P_\alpha^*\colon\alpha\in S\}\subset[\widehat\varrho\,]^{\varrho^+}$ 
is $\varrho$-a.d. But here $S\subset S_{\varrho^+}^{\widehat\varrho}$ 
is stationary and so by GCH we also have $\diamondsuit(S)$, so Lemma 2.1 
applies and hence $M(\widehat\varrho,\varrho^+,\varrho)\nrightarrow B$ 
holds in $V^P$.

\newpage

\centerline{{\bf \S 3. A ``stick''-like principle}} 
\smallskip

The aim of this section is to introduce a ``stick''-like combinatorial 
principle that will play an essential role in the proof of theorem 
1.4. We also look at some other results of purely
combinatorial nature and thus separate the combinatorial arguments from 
the rest, to be given in the next section.

\bigskip

{\bf 3.1. Definition.}\quad If $\kappa>\lambda\ge\mu\ge\omega$ then we 
denote by $\stick (\kappa,\lambda,\mu)$ the following statement:
There is a $\mu$-a.d. family $\irta\subset[\kappa]^\lambda$ with $|\irta|=
\kappa$ such that for every set $X\in[\kappa]^\kappa$ there is some $A
\in\irta$ with $A\subset X$; if $\irta$ is like this then we say that 
$\irta$ is a $\stick(\kappa,\lambda,\mu)$-family.

\bigskip

The relevance of this principle to our subject, in particular to 1.4, becomes
clear from the following result.

\bigskip

{\bf 3.2. Lemma.}\quad $\stick(\kappa,\lambda,\mu)$ implies that
$$
M(\kappa,\kappa,\mu)\nrightarrow B(\lambda).
$$

\bigskip

{\it Proof.}\quad Let $\irta$ be a $\stick(\kappa,\lambda,\mu)$-family and 
fix a partition $\{X_\xi\colon\xi\in\kappa\}\subset[\kappa]^\kappa$ of 
$\kappa$ into $\kappa$-many sets of size $\kappa$. Then we set 
$$
\irtb=\{A\in\irta\colon(\forall\xi\in\kappa)(|A\cap X_\xi|\le 1)\}.
$$
Clearly $|\irtb|=|\irta|=\kappa$, hence we may also fix a one-one 
enumeration $\irtb=\{B_\xi\colon\xi\in\kappa\}$ of $\irtb$.
Now, for every $\xi\in\kappa$ we set
$$
Y_\xi=X_\xi\cup B_\xi.
$$
Then it is obvious that the family 
$$
\irty=\{Y_\xi\colon\xi\in\kappa\}\subset[\kappa]^\kappa
$$
is $\mu$-a.d., hence we shall be done if we can show that $\irty$ has no
$\lambda$-transversal.

So assume that $T$ is such that 
$T\cap Y_\xi\ne\emptyset$ for all $\xi\in\kappa$. We claim that then the set 
$a=\{\xi\in\kappa\colon T\cap X_\xi\ne\emptyset\}$ has size $\kappa$.

Assume, indirectly, that $|a|<\kappa$. It is clear
that for any set $H\in [\kappa]^\kappa$, which satisfies 
$|H\cap X_\xi|\le 1$ for all $\xi\in\kappa$, we have
$$
|\{B_\xi\in\irtb\colon B_\xi\subset H\}|=\kappa.
$$
In particular, if $\alpha_\xi$ is the minimal member of $X_\xi$
for any $\xi\in\kappa$, then we may apply the above observation to the set
$$
H=\{\alpha_\xi\colon\xi\in\kappa\setminus a\}\in[\kappa]^\kappa
.$$
So there is some $\xi\in\kappa\setminus a$ such that $B_\xi\subset H$. 
But then, by the definition of the set $a$, we have both $T\cap 
H=\emptyset$, hence $T\cap B_\xi=\emptyset$ and $T\cap X_\xi=\emptyset$, 
i.e. $T\cap Y_\xi=\emptyset$, a contradiction.

But now, let us pick for every $\xi\in a$ an element $\beta_\xi\in T\cap 
X_\xi$ and set $K=\{\beta_\xi\colon\xi\in a\}$. We may then apply the above 
observation to the set $K\in[\kappa]^\kappa$ and find $B_\xi\in\irtb$ with 
$B_\xi\subset K$. So we conclude that $T\cap Y_\xi\supset B_\xi$, hence 
$|T\cap Y_\xi|\ge|B_\xi|=\lambda$, i.e. $T$ is not a $\lambda$-transversal.

\bigskip

{\it Remark.}\quad We have actually shown that $\irty$ has the following
stronger property: For any set $T$, if $|\{\xi\in\kappa\colon
 T\cap Y_\xi=\emptyset\}|<\kappa$ then there is some $Y_\xi\in\irty$ 
with $|T\cap Y_\xi|\ge\lambda$.

\bigskip

Our next result yields a method for ``stepping down'' in the second 
parameter $\lambda$ of a negative relation of the form $M(\kappa,\lambda, 
\mu)\nrightarrow B(\sigma)$. 

\bigskip

{\bf 3.3. Lemma.} Assume that $\tau<\lambda$ and we have both 

$(*)\qquad\qquad
\qquad M(\kappa,\lambda,\mu)\nrightarrow B(\sigma)$

and

$(**)\qquad\qquad\qquad M(\kappa,\lambda,\mu)\rightarrow B(\tau^+)$.

Then we also have
$$
M(\kappa,\tau,\mu)\nrightarrow B(\sigma).
$$

\bigskip

{\it Proof.}\quad Let $\irty=\{Y_\xi\colon\xi\in\kappa\}\subset
[\kappa]^\lambda$ be a $\mu$-a.d. family with no $\sigma$-transversal. 
With transfinite recursion on $\alpha\in\tau$ we define sets $T_\alpha$ 
that are all $\tau^+$-transversals of $\irty$ as follows.

Let $T_0$ be any $\tau^+$-transversal of $\irty$, it exists by ($**$). 
If $T_\beta$ has been defined for each $\beta\in\alpha\in\tau$ then 
for every $Y_\xi\in\irty$ we have 
$|Y_\xi\setminus \cup\{T_\beta\colon\beta\in\alpha\}|=\lambda$ 
because, by the inductive hypothesis, $|Y_\xi\cap T_\beta|\le\tau $ 
for each $\beta\in\alpha$. So we may now apply $(**)$ to the family
$\irty_\alpha=\{Y_\xi\setminus\cup\{T_\beta\colon\beta\in\alpha\}
\colon\xi\in\kappa\}$ and obtain a $\tau^+$-transversal $T_\alpha$ of 
$\irty_\alpha$ and hence of $\irty$.

Having completed the recursion, set
$$
T=\cup\{T_\alpha\colon\alpha\in\tau\}
$$
and $Z_\xi=Y_\xi\cap T$ for each $\xi\in\kappa$. It is clear from 
the construction that
$$
\irtz=\{Z_\xi\colon\xi\in\kappa\}
$$
is a $\mu$-a.d. subfamily of $[\kappa]^\tau$, so we'll be done if we can show
that $\irtz$ has no $\sigma$-transversal.

Since $\cup\irtz\subset T$, it suffices to show that if 
$U\subset T$ intersects every member of $\irtz$ then $|U\cap Z_\xi|\ge
\sigma$ for some $Z_\xi\in\irtz$. However, we know that there is a 
$\xi\in\kappa$ with $|U\cap Y_\xi|\ge \sigma$ which by $U\subset T$ and 
$Z_\xi=T\cap Y_\xi$ implies $|U\cap Z_\xi|\ge\sigma$, completing the proof.

Putting 1.2 and 3.3 together we immediately obtain the following result.

\bigskip

{\bf 3.4. Corollary.} (GCH)\qquad If
$$
M(\kappa,\kappa,\varrho)\nrightarrow B(\varrho^+)
$$
then for any $\lambda$ with $\varrho^+<\lambda<\kappa$ we have 
$$
M(\kappa,\lambda,\varrho)\nrightarrow B(\varrho^+)
$$
as well.

This implies that to prove 1.4 it suffices to concentrate on 
$M(\widehat\varrho,\widehat\varrho,\varrho)\nrightarrow B(\varrho^+)$, 
and so, by 3.2, on $\stick(\widehat\varrho,\varrho^+,\varrho)$.

Let us now make a few observations about the principles $\stick(\kappa, 
\lambda,\mu)$ that are less closely related to the main subject matter of this 
paper.

If $\stick(\kappa,\lambda,\mu)$ is valid then we obviously have a 
$\stick(\kappa,\lambda,\mu)$ family $\irta$ such that tp$A=\lambda$ 
for every $A\in\irta$. Let us now put
$$
S_\irta=\{\cup A\colon A\in\irta\},
$$ 
so $S_\irta\subset S_\varrho^\kappa$, where $\varrho=cf(\lambda) 
\le\lambda<\kappa$. We claim that if $\kappa$ is regular then 
$S_\irta$ is also stationary. Indeed, if $C\subset\kappa$ is c.u.b. then, 
as $|C|=\kappa$, there is some $A\in\irta$ with $A\subset C$ and thus 
$$
\cup A\in S_{\irta}\cap C\ne\emptyset.
$$
So, if GCH holds then we also have $\diamondsuit (S_{\irta})$, 
consequently from 2.1 and 1.2 we easily obtain the following result. 

\bigskip

{\bf 3.5. Proposition.} (GCH)\quad If $\kappa$ is regular then 
 $\stick(\kappa, \lambda,\mu)$ implies 
$M(\kappa,\lambda,\mu)\nrightarrow B$. 
Hence if $\kappa>\lambda>\varrho^+$ where $\kappa$ and $\varrho$ are 
regular then $\stick(\kappa,\lambda,\varrho)$ is false.

Thus, under GCH, for regular $\kappa$ 
and $\varrho$ the best we may hope for
is $\stick(\kappa,\varrho^+,\varrho)$, moreover, in view of 1.1, 
$\widehat\varrho$ is the smallest possible value for a $\kappa$ where 
this may happen. In fact, as follows from the next result,  
$\stick(\kappa,\varrho^+,\varrho)$ will fail for ``most'' regular 
$\kappa>\varrho^+$ even in ZFC.

\bigskip

{\bf 3.6. Proposition.}\quad If $\kappa$ is regular and for every 
$\lambda<\kappa$ we have $\lambda^\varrho<\kappa$ then 
$\stick(\kappa,\varrho^+,\varrho)$ is false.

\bigskip

{\it Proof.}\quad Assume that $\irta\subset[\kappa]^{\varrho^+}$ is 
$\varrho$-a.d. with tp$A=\varrho^+$ for all $A\in\irta$. According to 
what we have seen above, if we can show that $S_\irta$ is non-stationary 
in $\kappa$ then we are done.

Assume, indirectly, that $S_\irta$ is stationary and 
for each $\alpha\in S_\irta$ let $A_\alpha\in\irta$ be such that 
$\cup A_\alpha=\alpha$. For every $\alpha\in S_\irta$ let $f(\alpha)$ 
be the $\varrho^{\text{th}}$ element of $A_\alpha$, then 
$f$ is a regressive function on $S_\irta$ so by Fordor's 
theorem we have a stationary set $S\subset S_\irta$ and a $\gamma\in\kappa$ 
with $f(\alpha)=\gamma$ for every $\alpha\in S$. But then, using 
$|\gamma|^\varrho<\kappa$, we clearly have distict $\alpha,\beta\in S$ with 
$\gamma\cap A_\alpha=\gamma\cap A_\beta$, hence 
$|A_\alpha\cap A_\beta|\ge\varrho$, contradicting that $\irta$ is 
$\varrho$-a.d.

\bigskip

{\it Remark.}\quad The above argument actually yields the following stronger 
result: Under the assumptions of 3.6 even 
$\stick(\kappa,\varrho+\omega,\varrho)$ is false, with the obvious 
interpretion of this symbol.

\bigskip

Thus we have arrived ``down'' to $\stick(\kappa,\varrho,\varrho)$
that is ``easy'' to satisfy, being e.g. 
a consequence of the appropriate version of $\clubsuit$ at $\kappa$ and 
$\varrho$. In fact, in many cases it holds even in ZFC.

We close this section with two simple results concerning the behaviour of 
$\stick(\kappa,\lambda,\mu)$ in forcing extensions. The first one is a
preservation result.

\bigskip

{\bf 3.7. Proposition.}\quad Assume $\stick(\kappa,\lambda,\mu)$ where 
$\kappa$ is regular and $P$ is a forcing notion with $|P|<\kappa$ 
such that both $\lambda$ and $\mu$ remain cardinals in $V^P$ 
($\kappa$ does so automatically). Then $\stick(\kappa,\lambda,\mu)$ 
remains valid in $V^P$.

\bigskip

{\it Proof.}\quad Let $\irta$ be a $\stick(\kappa,\lambda,\mu)$-family 
in $V$. Now $|P|<\kappa=cf(\kappa)$ clearly implies that if 
$X\in[\kappa]^\kappa$ in $V^P$ then there is a $Y\in[X]^\kappa 
\cap V$, hence $A\subset Y\subset X$ for some $A\in\irta$, i.e. $\irta$ 
remains a $\stick(\kappa,\lambda,\mu)$-family in $V^P$.

The second result gives us a method to obtain $\stick(\chi,\varrho^+,
\varrho)$ for a given regular cardinal $\varrho$, assuming that we have
$\stick(\chi,\mu^+,\mu)$ for a singular cardinal $\mu$ of cofinality 
$\varrho$.

\bigskip

{\bf 3.8. Proposition.}\quad Assume $\stick(\chi,\mu^+,\mu)$, where 
$cf(\mu)=\varrho$, $\mu^{<\varrho}=\mu$, and $2^\mu<\chi=cf(\chi)$. 
Then $\stick(\chi,\varrho^+,\varrho)$ holds in $V^{Q(\mu)}$.

\bigskip

{\it Proof.}\quad Let $\irta\subset[\chi]^{\mu^+}$ be a 
$\stick(\chi,\mu^+,\mu)$-family in the ground model $V$. Then, in 
$V^{Q(\mu)}$, applying 2.3 we have for every $A\in\irta$ subset 
$A^*\in[A]^{\varrho^+}$ such that
$$
\irta^*=\{A^*\colon A\in\irta\}
$$
is $\varrho$-a.d. We claim that $\irta^*$ is a $\stick(\chi,
\varrho^+,\varrho)$-family. Since, by 2.2 (iii), we have 
$|Q(\mu)|\le 2^\mu<\chi$, similarly as in the proof of 3.7, every set 
$X\in[\chi]^\chi$ in $V^{Q(\mu)}$ has a ground model subset $Y$ 
with $|Y|=|X|=\chi$. But then there is an $A\in\irta$ with
$$
A^*\subset A\subset Y\subset X,
$$
and the proof is completed.

\bigskip\bigskip

\centerline{{\bf \S 4. The proof of 1.4}}
\smallskip

Assume GCH and that $\varrho=cf(\varrho)<\kappa$, where $\kappa$ is 
2-huge, in fact what we really need is the following property of $\kappa$ 
that is just a little more than being 1-huge:

There is an elementary embedding $j\colon V\to M$ with crit$(j)=\kappa$, 
$j(\kappa)=\lambda$ and $M^{\lambda^{(+\varrho+3)}}\subset M$, or 
equivalently there is a $\kappa$-complete normal ultrafilter 
$\irtd^*$ over $\irtp(H(\lambda^{(+\varrho+3)}))$ such that 
$$
\{M\colon M\prec H(\lambda^{(+\varrho+3)})\quad\&\quad 
M\cong H(\kappa^{(+\varrho+3)})\}\in\irtd^*.
$$
We shall be working with the projection $\irtd$ of $\irtd^*$ to 
$ H(\lambda^{(+\varrho+1)})$, i.e.
$$
\irtd=\{A\subset H(\lambda^{(+\varrho+1)})\colon
\{a\subset H(\lambda^{(+\varrho+3)})\colon
a\cap H(\lambda^{(+\varrho+1)})\in A\}\in\irtd^*\}.
$$
Then, of course, $\irtd$ is a $\kappa$ complete normal ultrafilter 
over $\irtp( H(\lambda^{+\varrho+1}))$ such that 
$$
X=\{M\colon M\prec H(\lambda^{(+\varrho+1)})\quad\&\quad  
M\cong H(\kappa^{(+\varrho+1)})\}\in\irtd.
$$
Let us write, for simplicity, $\kappa^{(+\varrho)}=\mu$ and 
$\lambda^{(+\varrho+1)}=\chi$. Combining the above with Solovay's result
as in the final part of section 2, we conclude that there is a 
stationary set $S\subset S_{\mu^+}^\chi$ such that 
for each $\delta\in S$ we have 
$M_\delta\in X$, $\cup(M_\delta\cap\chi)=\delta$, moreover 
$\{M_\delta\colon\delta\in S\}\in\irtd$ is $\mu$-a.d.. In what folows, 
we shall write $Y_\delta=M_\delta\cap\chi$ for $\delta\in S$.

The crucial part of our proof is the following result.

\bigskip

{\bf 4.1. Lemma.}\quad There is a sequence 
$\langle f_\delta\colon\delta\in S\rangle$ such that
\begin{itemize}
\item[(i)] $f_\delta\colon Y_\delta\to Y_\delta$ for each 
$\delta\in S$;
\item[(ii)] for every $f\colon\chi\to\chi$ the set
$$
\{\delta\in S\colon f_\delta\subset f\}
$$
\end{itemize}
is stationary in $\chi$.

\bigskip

{\it Proof.}\quad Let us write, for $\delta\in S$,
$$
Y_\delta=\{\alpha_{\delta,\zeta}\colon\zeta\in\mu^+\},
$$
the increasing enumeration of $Y_\delta$. The functions 
$f_\delta\colon Y_\delta\to Y_\delta$ will be defined by a simple 
transfinite recursion in such a way that for each $\delta\in S$ the set 
$H_\delta=\{\zeta\in\mu^+\colon f_{\alpha_{\delta,\zeta}}\restriction
 Y_\delta\cap Y_{\alpha_{\delta,\zeta}}\subset f_\delta\}$
be non-stationary in $\mu^+$, if this is possible at all.

All we have to do now is to check that (ii) holds.  Assume, indirectly, 
that $f\colon\chi\to\chi$ and $C\subset\chi$ c.u.b. exist such that 
$f_\delta\not\subset f$ for every $\delta\in S\cap C$.

For any $\alpha\in\chi$, as normality of $\irtd$ implies its fineness, 
we have
$$
A_\alpha=
\{M_\delta\colon
\{\alpha,f_\delta(\alpha)\}\subset M_\delta\}\in\irtd.
$$
Also, for any pair $\langle\alpha,\beta\rangle\in\chi^2$ we can define 
$A_{\alpha,\beta}\in\irtd$ so that
$$A_{\alpha,\beta}=\cases
\{M_\delta\colon f_\delta(\alpha)=\beta\} &\\
\text{or}&\\
\{M_\delta\colon f_\delta(\alpha)\ne\beta\}. &
\endcases
$$
Then, by the normality of $\irtd$, there is a (clearly stationary) 
subset $S_1\subset S\cap C$
such that 
$$
X_1=\{M_\delta\colon\delta\in S_1\}\in\irtd
$$
and if $\delta\in S_1$, $\langle\alpha,\beta\rangle\in Y_{\delta}^2$ then 
$M_\delta\in A_\alpha\cap A_{\alpha,\beta}$.

Let $M_\delta\in X_1\cap A_\alpha\in\irtd$ where $\alpha\in\chi$, 
then clearly $g(\alpha)=f_\delta(\alpha)$ does not depend on 
$\delta$, moreover
$$
\{M_\delta\colon f_\delta(\alpha)=g(\alpha)\}\in\irtd.
$$
This implies that for every $\alpha\in\chi$ we have
$$
A_{\alpha,g(\alpha)}=\{M_\delta\colon f_\delta(\alpha)=g(\alpha)\},
$$
consequently $f_\delta\subset g$ whenever $\delta\in S_1$. In particular, 
as $S_1\subset C$ and $f_\delta\not\subset f$ for $\delta\in C$, we have 
$f\ne g$.

Now, applying the normality of our original ultrafilter $\irtd^*$, we 
can find $N\prec H(\chi^{++})$ such that $N\cong H(\mu^{+++})$, moreover 
\begin{itemize}
\item[(a)] $\langle M_\delta\colon\delta\in S\rangle$,
$\langle f_\delta\colon\delta\in S\rangle$, $S_1,f,g,C,\irtd\in N$;
\item[(b)] for any $Z\in N\cap\irtd$ we have $N\cap H(\chi)\in Z\cap X_1$.
\end{itemize}

Let $h\colon N\to H(\mu^{+++})$ be the Mostowski collapse, then 
$h(\chi)=\mu^+$. Moreover, from (a) and (b) it follows that 
$N\cap H(\chi)=M_{\delta^*}$, where 
 $\delta^*\in S_1$. By elementarity $N\models$ 
``$S_1$ is stationary in $\chi$'', hence $h(S_1)$ is stationary in 
$h(\chi)=\mu^+$, or in the other words the set
$$
H=\{\zeta\in\mu^+\colon\alpha_{\delta^*,\zeta}\in S_1\}
$$
is stationary in $\mu^+$. But if $\alpha_{\delta^*,\zeta}\in S_1$ 
then we have $f_{\alpha_{\delta^*,\zeta}}\subset g$ as well as 
$f_{\delta^*}\subset g$, hence $f_{\alpha_{\delta^*,\zeta}}
\restriction Y_{\delta^*}\cap Y_{\alpha_{\delta^*,\zeta}}
\subset f_{\delta^*}$. So we conclude from $H\subset H_{\delta^*}$ that 
at step $\delta^*$ of the transfinite construction we could not make 
$H_{\delta^*}$ non-stationary.

However, as $f,g,C\in N$ we have on one hand that 
$f_N=f\restriction Y_{\delta^*}\colon Y_{\delta^*}\to Y_{\delta^*}$, 
moreover the set
$$
\{\zeta\in\mu^*\colon\alpha_{\delta^*,\zeta}\in C\}=h(C)
$$
is c.u.b. in $\mu^+$. By elementarity, as $f\ne g$, for every $\zeta\in h(C)$ 
there is a $\gamma\in N\cap M_{\alpha_{\delta^*,\zeta}}$ such that 
$f(\gamma)\ne g(\gamma)=f_{\alpha_{\delta^*,\zeta}}(\gamma)$, i.e.
$$
f_{\alpha_{\delta^*,\zeta}}\restriction Y_{\delta^*}\cap
 Y_{\alpha_{\delta^*,\zeta}}\not\subset f_N.
$$
This, however contradicts our above conclusion because $f_N$ would make, 
at step $\delta^*$, the set $H_{\delta^*}$ non-stationary in $\mu^*$.

Now from 4.1 we easily obtain the following result, where the notation 
is the same. 

\bigskip

{\bf 4.2. Proposition.}\quad $\stick(\chi,\mu^+,\mu)$ is valid. 

\bigskip

{\it Proof.}\quad Let $S^*=\{\delta\in S\colon f_\delta
\text{ is strictly increasing}\}$ and for each $\delta\in S^*$ let 
$Z_\delta=f_\delta{}''Y_\delta$. We claim that 
$\irtz=\{Z_\delta\colon\delta\in S^*\}\subset [\chi]^{\mu^+}$ is a 
$\stick(\chi,\mu^+,\mu)$-family. Since $Z_\delta\subset M_\delta$, 
$\irtz$ is clearly $\mu$-a.d. 
Now, for any set $Z\in[\chi]^\chi$ let $f$ be its increasing 
enumerating function, then
$$
S_f=\{\delta\colon f_\delta\subset f\}
$$
is stationary and also $S_f\subset S^*$. But for any $\delta\in S_f$ 
we clearly have $Z_\delta\subset Z$.

Now, it is very easy to complete the proof of 1.4. First note that 3.8 may be 
applied, i.e. in $V^{Q(\mu)}$ we have $\stick (\chi,
\varrho^+,\varrho)$. Next, similarly as in \S 2, if one collapses $\lambda$ 
to $\varrho^{++}$ in $V^{Q(\mu)}$ using $L v(\lambda,\varrho^{++})$ 
then the forcing $P=Q(\mu)* Lv(\lambda,\varrho^{++})$ 
is as required because it is $\varrho$-complete, preserves GCH, moreover 
$\stick(\widehat\varrho,\varrho^+,\varrho)$ holds true in $V^P$. Indeed, 
the last part follows because $\chi=\widehat\varrho$ in $V^P$ and 
$\stick (\chi,\varrho^+,\varrho)$ is preserved by the Levy-collapse, 
using 3.7 and $|Lv(\lambda,\varrho^{++})|<\chi$.

\newpage

\centerline{{\bf References}}

\bigskip

\begin{itemize}
\item[[BDJShSz]] Z. T. Balogh, S. W. Davis, W. Just, S. Shelah and 
J. Szept\'ycki, Strongly almost disjoint sets and weakly uniform bases, 
Preprint no. 12 (1997/98), Hebrew Univ. Jerusalem, Inst. of Math.
\smallskip
\item[[EH]] P. Erd\H os and A. Hajnal, On a property of families of sets, 
Acta Math. Acad. Sci. Hung., {\bf 12}(1961), p.87--124.
\smallskip
\item[[HJSh]] A. Hajnal, I. Juh\'asz and S. Shelah, Splitting strongly 
almost disjoint families, TAMS {\bf 295}(1986), p. 369--387.
\smallskip
\item[[S]] R.. Solovay, Strongly compact cardinals and the GCH, 
Proc. Symp. in Pure Math., vol.XXV,  1974, p. 365--372.
\smallskip
\item[[W]] N. H. Williams, Combinatorial Set Theory, Studies in Logic, 
Vol. 91, North-Holland, Amsterdam, 1977.  
\end{itemize}

\vskip2truecm
\baselineskip=12pt
\noindent Authors' addresses:\newline
\bigskip

\noindent 1) Department of Mathematics\newline
Rutgers University\newline
New Brunswick, New Jersey 08903\newline
U.S.A.\newline
E-mail: ahajnal@math.rutgers.edu
\bigskip

\noindent 2) Mathematical Institute\newline
of the Hungarian Academy of Sciences\newline
P.O.Box 127\newline
1364 Budapest, Hungary\newline
E-mail: juhasz@math-inst.hu
\bigskip

\noindent 3) Department of Mathematics\newline
Rutgers University\newline
New Brunswick, New Jersey 08903\newline
U.S.A.\newline
and\newline
Institute of Mathematics\newline
The Hebrew University\newline
91904 Jerusalem, Israel\newline
E-mail: shelah@math.huji.ac.il

\end{document}